\documentclass[11pt]{amsart}
\usepackage{amssymb,amsmath,amsfonts,amscd,euscript,mathrsfs}
\usepackage[all]{xy}
\usepackage{color}
\usepackage{tikz-cd}

\usepackage[backref=page]{hyperref}
\usepackage{todonotes}

\usepackage{amsthm,enumitem,bm,pict2e,xcolor}
\usepackage{indentfirst}

\newcommand{\nc}{\newcommand}

\numberwithin{equation}{section}
\newtheorem{thm}{Theorem}[section]
\newtheorem{prop}[thm]{Proposition}
\newtheorem{lem}[thm]{Lemma}
\newtheorem{cor}[thm]{Corollary}
\newtheorem{rem}[thm]{Remark}

\newtheorem{example}[thm]{Example}

\newtheorem{theoremA}{Theorem}

\nc{\la}{\lambda}
\nc{\al}{\alpha }
\nc{\be}{\beta }
\nc{\om}{\omega }

\nc{\ch}{{\mathop {\rm ch}}}
\nc{\ld}{\ldots}
\nc{\cd}{\cdots}
\nc{\hk}{\hookrightarrow}
\nc{\T}{\otimes}
\nc{\Gr}{\mathrm{Gr}}
\nc{\Fl}{\mathrm{Fl}}
\nc{\Spec}{\mathrm{Spec}}
\nc{\Proj}{\mathrm{Proj}}
\nc{\ov}{\overline}
\nc{\msl}{\mathfrak{sl}}
\nc{\mgl}{\mathfrak{gl}}
\nc{\U}{\mathrm U}

\newcommand{\bC}{{\mathbb C}}
\newcommand{\bR}{{\mathbb R}}
\newcommand{\bZ}{{\mathbb Z}}

\newcommand{\bP}{{\mathbb P}}
\newcommand{\bA}{{\mathbb A}}

\newcommand{\eO}{\EuScript{O}}

\newcommand{\bfx}{{\bf x}}
\newcommand{\bfy}{{\bf y}}
\newcommand{\bfs}{{\bf s}}

\newcommand{\SSYT}{\mathrm{SSYT}}

\begin{document}

\title[Specializations of frozen variables]
{Specializations of frozen variables for Grassmannian cluster algebra}

\author{Evgeny Feigin}
\address{Evgeny Feigin:\newline
School of Mathematical Sciences, Tel Aviv University, Tel Aviv, 69978, Israel}
\email{evgfeig@gmail.com}

\author{Jian-Rong Li}
\address{Jian-Rong Li:\newline 
Faculty of Mathematics, University of Seville, Calle Tarfia s/n, 41012 Seville, Spain}
\email{lijr07@gmail.com}

\begin{abstract}
The homogeneous coordinate ring of a Grassmannian variety $\mathrm{Gr}_{k,n}$ is equipped with a
structure of a cluster
algebra. The frozen variables are identified with the Pl\"ucker coordinates corresponding
to length-$k$ segments of ${\mathbb Z}_n$. We study two families of algebras: the first family
is obtained by specializing all the frozen variables and the second corresponds to the 
specialization of the nonwrapping variables. We construct flat bases in both families and
study their affine spectra. The special fibers of the families are cones over projective varieties. 
We identify the irreducible components of these projective varieties with certain positroid varieties.
\end{abstract}

\maketitle

\section{Introduction}

Let $A$ be a cluster algebra with frozen variables $f_1,\dots,f_n$. In particular, 
$A$ is an algebra over the polynomial ring $S$ in variables $f_i$ (we do not invert
the frozen variables). Hence, one gets the embedding $S\to A$ and the projection
$\Spec(A)\to \bA^n=\Spec(S)$ -- the frozen variables specialization map. 
In particular, the preimage of zero can be understood as the deepest strata of the
corresponding cluster variety (see \cite{BF19,CMNC22,GHKK18}).
In this paper we are interested in the fibers of the frozen variables specialization
map for the Grassmannian cluster algebra. 

Let $R_{k,n}$ be the Grassmann algebra. It can be realized as the homogeneous
coordinate ring of the Grassmannian $\Gr_{k,n}$ with respect to the classical 
Pl\"ucker embedding, i.e. $\Proj(R_{k,n})\simeq \Gr_{k,n}$. Similarly, one 
can take the affine spectrum $\Spec(R_{k,n})$ and arrive at the affine cone 
$\widehat\Gr_{k,n}$. Explicitly, the Grassmann algebra is
isomorphic to a quotient of the polynomial ring in Pl\"ucker variables 
$P_I$ ($I$ running over the cardinality $k$ subsets of the set $\{1,\dots,n\}$) modulo
a quadratic ideal. 

It was shown 
in \cite{Sc06} that $R_{k,n}$ admits a structure of a cluster algebra. In particular,
variables $P_{C_a}$ for $C_a=\{a,\dots,a+k-1\}$, $1\le a\le n$ are the frozen variables (the indices
$a+i$ are understood modulo $n$). Let us denote these Pl\"ucker variables by $f_a$
and let $S$ be the polynomial ring in variables $f_a$. 
We obtain a map $\widehat\Gr_{k,n}\to \bA^n=\Spec(S)$; the preimage over a point 
$\bfx = (x_1,\dots,x_n)\in \bA^n$ is equal to $\Spec(R_{k,n}^\bfx)$, where
\begin{equation}\label{eq:intrcycfam}
R_{k,n}^\bfx = R_{k,n}/(f_1-x_1,\dots,f_n-x_n).
\end{equation}

The family can be seen as a degeneration of a generic fiber to the 
most special one corresponding to  $x_i=0$ for all $i$. The generic case
has to do with the representation theory of the quantum affine algebras.  
The most special
case $\bfx=0$ is the most interesting one and is in some sense responsible for 
the behavior of the whole family.

In order to study this degeneration, it is useful to consider another degeneration direction.
Namely, the algebras $R_{k,n}$ admit Gelfand--Tsetlin toric degeneration 
(see \cite{GL96,KM05,Stu93,SZ93}). The resulting algebra $\overline{R}_{k,n}$ is a quotient 
of the polynomial algebra in variables $\ov{P}_I$ by a binomial ideal. We denote the
analogues of the frozen variables, $\ov P_{a,\dots,a+k-1}$, by $\ov f_a$ (as before, 
$a=1,\dots,n$). One gets the toric analogue of the family \eqref{eq:intrcycfam}:
$\ov R_{k,n}^\bfx = \ov R_{k,n}/(\ov f_1-x_1,\dots,\ov f_n-x_n)$.
The polynomial algebra in variables $\ov f_a$ is denoted by $\ov S$.

We prove the following theorem. 
\begin{theoremA}
The polynomial algebra $S$ (resp., $\ov S$) acts freely
on $R_{k,n}$ (resp., on $\ov R_{k,n}$).
For all $\bfx\in\bA^n$, the following are flat bases
of $R^\bfx_{k,n}$: noncrossing monomials without
factors $P_{C_a}$, $a\in[n]$, cyclically reduced
semistandard monomials, and cyclically reduced
$\ch$-elements.
The algebra $\ov R^\bfx_{k,n}$ admits the first
two corresponding flat bases.
\end{theoremA}

There is a version of the family \eqref{eq:intrcycfam} which has a deep meaning in the 
representation theory of quantum affine algebras \cite{HL10,CDFL20}. More precisely, 
if $\ell=n-k-1$, Hernandez and
Leclerc introduced a monoidal subcategory $\mathcal C_\ell$ of
finite-dimensional representations of
$\U_q(\widehat{\mathfrak{sl}}_k)$, and there is an isomorphism
\[
\bC\otimes_{\bZ}K_0(\mathcal C_\ell)
\simeq
R_{k,n}/(f_1-1,\ldots,f_{n-k+1}-1).
\]
Under this correspondence, simple objects of $\mathcal C_\ell$ are
related to the dual canonical basis of
$R_{k,n}/(f_1-1,\ldots,f_{n-k+1}-1)$, and this correspondence was
made explicit in \cite{CDFL20} using semistandard tableaux. Motivated by this, 
we consider the polynomial ring ${\mathsf S}$ in variables
$f_1,\ldots,f_{n-k+1}$, keeping only those frozen variables
$f_a=P_{C_a}$ for which the consecutive indices
$a,\ldots,a+k-1$ do not exceed $n$. 
For $\bfy = (y_1,\dots, y_{n-k+1})\in \bA^{n-k+1}$, let 
\[
{\mathsf R}^\bfy_{k,n} = R_{k,n}/(f_1-y_1,\dots,f_{n-k+1}-y_{n-k+1}).  
\]
Similarly, one gets the corresponding GT version $\ov {\mathsf R}^\bfy_{k,n}$. Here, in order to distinguish these algebras from the cyclic families
$R_{k,n}^{\bfx}$ and $\ov R_{k,n}^{\bfx}$, we use the notation
${\mathsf R}_{k,n}^{\bfy}$ and $\ov{\mathsf R}_{k,n}^{\bfy}$ instead of
$R_{k,n}^{\bfy}$ and $\ov R_{k,n}^{\bfy}$.

\begin{theoremA}
The polynomial algebra ${\mathsf S}$ (resp., $\ov{\mathsf S}$)
acts freely on $R_{k,n}$ (resp., on $\ov R_{k,n}$).
For all $\bfy\in\bA^{n-k+1}$, the following are flat bases of
${\mathsf R}_{k,n}^{\bfy}$: noncrossing monomials without
factors $P_{C_a}$, $a\in[n-k+1]$, reduced semistandard
monomials, and reduced $\ch$-elements.
The algebra $\ov{\mathsf R}_{k,n}^{\bfy}$ admits the first two
corresponding flat bases.
\end{theoremA}

The following schematic diagrams illustrate the relations
between the algebras considered in this paper.
The horizontal arrows indicate specialization of the
frozen parameters to zero, while the vertical arrows
indicate passage to the corresponding Gelfand--Tsetlin
algebras. The arrows are not meant to represent
algebra homomorphisms between the fibers. The upper left corner represents a general fiber of the family
($\bfx\in\bA^n$, $\bfy\in\bA^{n-k+1}$).
\[
\begin{tikzcd}
R^\bfx_{k,n} \arrow[r] \arrow[d, "\mathrm{GT}"'] & R^0_{k,n} \arrow[d, "\mathrm{GT}"',swap] \\
\ov R^\bfx_{k,n} \arrow[r]           & \ov R^0_{k,n}
\end{tikzcd}\qquad
\begin{tikzcd}
{\mathsf R}^\bfy_{k,n} \arrow[r] \arrow[d, "\mathrm{GT}"'] & {\mathsf R}^0_{k,n} \arrow[d, "\mathrm{GT}"',swap] \\
\ov {\mathsf R}^\bfy_{k,n} \arrow[r]           & \ov {\mathsf R}^0_{k,n}
\end{tikzcd}
\]

We now turn to the geometric part of the story. The families of algebras $R_{k,n}^\bfx$ and
${\mathsf R}_{k,n}^\bfy$ define families of the corresponding affine algebraic varieties. The most 
interesting members of the family are the special fibers over zero, since they are the cones over 
certain projective varieties. In order to identify these underlying projective algebraic schemes, we use
the open positroid varieties of \cite{KLS13,KLS14}. The positroid varieties form a stratification of
$\Gr_{k,n}$. They are labeled by various combinatorial objects; in particular, they can be described 
via matroids of the nonvanishing Pl\"ucker coordinates.

\begin{theoremA}
The family $\Spec(R^\bfx_{k,n})$ is flat over $\bA^n$. 
For all $\bfx\in \bA^n$, the variety $\Spec(R^\bfx_{k,n})$ is connected, reduced, and of pure dimension 
$k(n-k)+1-n$. 
The generic fiber is smooth, rational, and irreducible of dimension given by $k(n-k)+1-n$. The special fiber $\Spec(R^0_{k,n})$ is a cone over 
an equidimensional projective variety $Y_{k,n}$. Each irreducible
component of $Y_{k,n}$ is a positroid variety. 
\end{theoremA}

Now let us consider the noncyclic case, i.e. the family ${\mathsf R}^\bfy_{k,n}$ over $\bA^{n-k+1}$.
It turned out that the varieties responsible for the geometry of the special fiber are the
interval positroid varieties \cite{Kn14,KL14,Paw18}. 

\begin{theoremA}
The family $\Spec({\mathsf R}^\bfy_{k,n})$ is flat over $\bA^{n-k+1}$.
For all $\bfy\in \bA^{n-k+1}$, the variety $\Spec({\mathsf R}^\bfy_{k,n})$ is connected, reduced, and of pure dimension 
$(k-1)(n-~k)$. 
If $y_j\ne 0$ for all $j$, then $\Spec({\mathsf R}^\bfy_{k,n})\simeq \bA^{(k-1)(n-k)}$. The special fiber $\Spec({\mathsf R}^0_{k,n})$ is a cone over 
an equidimensional projective variety ${\mathsf Y}_{k,n}$. The irreducible 
components of ${\mathsf Y}_{k,n}$ are the interval positroid varieties of codimension $n-k+1$ contained in $\{f_1=\dots=f_{n-k+1}=0\}$. 
The number of irreducible components is
$F_{n-k+1}^{(k)}-1$, where $F_m^{(k)}$
denotes the $m$-th $k$-step Fibonacci number. 
\end{theoremA}

Our paper is organized as follows. In Section \ref{sec:prelim} we collect the preliminary
material used in the main body of the paper. In Section \ref{sec:algcyc} we study the
family of algebras $R^\bfx_{k,n}$. Section \ref{sec:algnoncyc} is devoted to the study
of the family ${\mathsf R}^\bfy_{k,n}$. In Sections \ref{sec:varcyc} and \ref{sec:varnoncyc}
we describe geometric properties of the spectra of algebras $R^\bfx_{k,n}$
and ${\mathsf R}^\bfy_{k,n}$, including the special fibers. 

\section*{Acknowledgments}
EF was partially supported by the ISF grant 493/24. JL was supported by the Beatriz Galindo Senior grant BG24/00114 from the Spanish Ministry of Science, Innovation and Universities. 

\subsection*{AI use disclosure} ChatGPT 5.6 was used as an auxiliary tool to assist with computations, locate relevant references, and discuss bases of various quotients of Grassmannian cluster algebras. The manuscript was written by the authors, and all mathematical ideas presented in this paper are their own. The authors take full responsibility for the validity of the results and the correctness of the arguments.

\section{Preliminaries}\label{sec:prelim}

In this section, we recall the basic notions and results needed throughout the paper. We begin with the Grassmannian coordinate ring and its cluster algebra structure, and then recall the connection with quantum affine algebras and the dual canonical basis. We also review the Gelfand--Tsetlin toric degeneration, positroid varieties, and interval positroid varieties. 

We fix the base field $\bC$ and assume $2\le k\le n-2$ (the cases $k=1$ and $k=n-1$ are trivial, because all the Pl\"ucker variables are frozen, see subsection \ref{subsec:cluster-algebras}) .
We use the notation $[n]$ for the set $\{1,\dots,n\}$.
\subsection{Grassmannians}
\label{subsec:grassmannians}

Let $\Gr_{k,n}$ be the Grassmannian of $k$-dimensional subspaces of
$\bC^n$.  Under the Pl\"ucker embedding, $\Gr_{k,n}\hk \bP\!\left(\bigwedge^k\bC^n\right)$,
its homogeneous coordinate ring is the Grassmann algebra
$R_{k,n}:=\bC[\Gr_{k,n}]$.
For a $k$-subset $I=\{i_1<\cdots<i_k\}\subset [n]$,
we denote the corresponding Pl\"ucker coordinate by $P_I$.  Thus
$R_{k,n}$ is generated by the $P_I$'s subject to the quadratic
Pl\"ucker relations.  We write
\[
   \widehat{\Gr}_{k,n}:=\Spec(R_{k,n})
\]
for the affine cone over the Grassmannian.

We shall use the standard monomial basis of $R_{k,n}$.  A
semistandard tableau of rectangular shape with $k$ rows will be
identified with an ordered collection of columns
$T=(I_1,\ldots,I_d)$,
$I_j\in\binom{[n]}{k}$,
weakly increasing in each row of $T$.  Let $P_T:=P_{I_1}\cdots P_{I_d}$, then
monomials $P_T$ indexed by semistandard tableaux with $k$ rows form a $\bC$-basis of $R_{k,n}$; see, for example,
\cite{Hod43, LR08}.

For $a\in\bZ_n$, let
\[
   C_a=\{a,a+1,\ldots,a+k-1\}\subset\bZ_n,
\]
where the elements are written increasingly as a $k$-subset of
$[n]$.  The Pl\"ucker coordinates $P_{C_a}$ will play a distinguished
role throughout the paper.

\subsection{Cluster algebras}
\label{subsec:cluster-algebras}

We briefly recall the definition of a cluster algebra of geometric
type; see \cite{FZ02,FZ07} for the general theory and \cite{BFZ05}
for upper cluster algebras and related geometric realizations.
Let
\begin{align*}
\mathcal F=\bC(x_1,\ldots,x_r,f_1,\ldots,f_m).
\end{align*}
A seed is a pair
$\Sigma=(\widetilde{\bf x},\widetilde B)$, where
$\widetilde{\bf x}=(x_1,\ldots,x_r,f_1,\ldots,f_m)$,
the variables $x_1,\ldots,x_r$ are mutable,
$f_1,\ldots,f_m$ are frozen, and
$\widetilde B=(b_{ij})$ is an $(r+m)\times r$ integer matrix whose
principal $r\times r$ part is skew-symmetrizable.

For $1\leq j\leq r$, mutation in direction $j$ replaces
$\Sigma$ by
$\mu_j(\Sigma)=(\widetilde{\bf x}',\widetilde B')$, where
$\widetilde{\bf x}'$ is obtained from $\widetilde{\bf x}$ by
replacing $x_j$ with $x_j'$, leaving all other mutable and frozen
variables unchanged, and
\begin{equation}\label{eq:cluster-exchange-relation}
x_jx_j'
=
\prod_{i=1}^{r+m}x_i^{[b_{ij}]_+}
+
\prod_{i=1}^{r+m}x_i^{[-b_{ij}]_+},
\end{equation}
with $x_{r+a}=f_a$ and $[c]_+=\max(c,0)$.  The mutated matrix is
given by
\[
b'_{i\ell}
=
\begin{cases}
-b_{i\ell}, & i=j\text{ or }\ell=j,\\
b_{i\ell}+[b_{ij}]_+[b_{j\ell}]_+
-[-b_{ij}]_+[-b_{j\ell}]_+, & \text{otherwise}.
\end{cases}
\]
The mutable variables appearing in seeds obtained by repeated
mutations are the cluster variables.

We use the convention that frozen variables are not inverted.
Thus the cluster algebra $A=A(\Sigma)$ is the
$\bC[f_1,\ldots,f_m]$-subalgebra of $\mathcal F$ generated by
all cluster variables together with $f_1,\ldots,f_m$.
Set $S=\bC[f_1,\ldots,f_m]$.  The inclusion $S\hk A$ induces
\[
\Spec(A)\longrightarrow\Spec(S)=\bA^m,
\]
and the fiber over
${\bf a}=(a_1,\ldots,a_m)$ has coordinate ring
$$A/(f_1-a_1,\ldots,f_m-a_m).$$  In particular, the fiber over the
origin is $A/(f_1,\ldots,f_m)$.

Scott showed in \cite{Sc06} that $R_{k,n}$ carries a cluster algebra
structure of geometric type.   The frozen variables are
$f_a=P_{C_a}$, $a\in\bZ_n$ (recall that the entries of $C_a$ are written in the increasing order).  We write
\begin{equation}\label{eq:prelim-cyclic-family}
R_{k,n}^{\bfx}
=
R_{k,n}/(f_1-x_1,\ldots,f_n-x_n),
\qquad
\bfx=(x_1,\ldots,x_n)\in\bA^n.
\end{equation}

We also consider the noncyclic family, in which we retain only the nonwrapping consecutive
Pl\"ucker coordinates
$f_a=P_{a,a+1,\ldots,a+k-1}$,
$1\leq a\leq n-k+1$.  For
$\bfy=(y_1,\ldots,y_{n-k+1})\in\bA^{n-k+1}$, set
\begin{equation}\label{eq:prelim-noncyclic-family}
{\mathsf R}_{k,n}^{\bfy}
=
R_{k,n}/(f_1-y_1,\ldots,f_{n-k+1}-y_{n-k+1}).
\end{equation}
We use ${\mathsf R}_{k,n}^{\bfy}$ rather than
$R_{k,n}^{\bfy}$ in order to distinguish the noncyclic family from
the cyclic family $R_{k,n}^{\bfx}$.

\subsection{Quantum affine algebras and dual canonical basis}
\label{subsec:quantum-affine-algebras}

We recall results on quantum affine algebras and Grassmannians \cite{HL10, CDFL20} which motivate the noncyclic
specialization.  Set $\ell=n-k-1$, let
$I=[1,k-1]$, and consider the quantum affine algebra
$\U_q(\widehat{\msl}_k)$.  Let $\mathcal P^+$ be the free
commutative monoid generated by the $\ell$-weight variables
$Y_{i,s}$, $i\in I$, $s\in\bZ$.  We use the height function
$\xi:I\to\bZ$ given by $\xi(i)=i-2$, and let
$\mathcal P_\ell^+$ be the submonoid of $\mathcal P^+$ generated
by $Y_{i,\xi(i)-2r}$, $i\in I$, $0\leq r\leq\ell$. Let $\mathcal C_\ell$ be the full subcategory of
finite-dimensional $\U_q(\widehat{\msl}_k)$-modules
whose simple composition factors have highest
$\ell$-weights in $\mathcal P_\ell^+$. After extension of
scalars to $\bC$, the Grothendieck ring $K_0(\mathcal{C}_{\ell})$ of $\mathcal C_\ell$ with
the Grassmannian quotient
\[
\bC\otimes_{\bZ}K_0(\mathcal C_\ell)
\ \simeq\
{\mathsf R}_{k,n}^{\mathbf 1}
=
R_{k,n}/(f_1-1,\ldots,f_{n-k+1}-1).
\] 

We next recall from \cite{CDFL20} the tableau parametrization of the
dual canonical basis of ${\mathsf R}_{k,n}^{\mathbf 1}$.
Let ${\SSYT}(k,[n])$ denote the set of semistandard rectangular tableaux with
$k$ rows and entries in $[n]$.  It is a commutative cancellative
monoid under row-wise union, denoted by $T\cup T'$.  A tableau is
called \emph{trivial} if it is a product of the nonwrapping consecutive
columns $C_a=(a,a+1,\ldots,a+k-1)$, $1\leq a\leq n-k+1$.
Two tableaux $T,T'$ are equivalent if there exist trivial tableaux
$U,U'$ such that $T\cup U=T'\cup U'$.  We denote the resulting
quotient monoid by ${\SSYT}(k,[n],\sim)$.  Each equivalence class
contains a unique reduced representative, namely one having no
nontrivial trivial factor.

By \cite[Theorem~3.17]{CDFL20}, there is a monoid
isomorphism
$\mathcal P_\ell^+
   \xrightarrow{\sim}
   \mathrm{SSYT}(k,[n],\sim)$. 
Hence simple objects of $\mathcal C_\ell$ may be parametrized by
reduced semistandard tableaux.

For $\overline T\in\mathrm{SSYT}(k,[n],\sim)$, denote by
$\widetilde{\ch}(\overline T)\in{\mathsf R}_{k,n}^{\mathbf 1}$
the element corresponding to the associated simple module.  Then
\[
\left\{\widetilde{\ch}(\overline T):
\overline T\in\mathrm{SSYT}(k,[n],\sim)\right\}
\]
is the dual canonical basis of ${\mathsf R}_{k,n}^{\mathbf 1}$.
More explicitly, let $T'\sim T$ be the small-gap representative of
$\overline T$, let $d$ be the number of columns of $T'$, and let
$w_T\in S_d$ be as in \cite[Section~5]{CDFL20}.  Then
\cite[Theorem~5.8]{CDFL20} gives
\begin{equation}\label{eq:prelim-CDFL-KL-formula}
\widetilde{\ch}(\overline T)
=
\sum_{u\in S_d}
(-1)^{\ell(uw_T)}
p_{uw_0,w_Tw_0}(1)\,P_{u;T'},
\end{equation}
where $p_{v,w}$ denotes the Kazhdan--Lusztig polynomial and
$P_{u;T'}$ is the Pl\"ucker monomial defined in
\cite[Definition~5.6]{CDFL20}.

In \cite[Definition~5.9]{CDFL20}, a homogeneous Laurent lift
is associated with each $T\in{\SSYT}(k,[n])$ as follows.
By \cite[Lemma~3.13 and Remark~3.14]{CDFL20},
there is a unique small-gap tableau $T'$
equivalent to $T$, and there exist trivial tableaux
$A$ and $B$ such that $T\cup A=T' \cup B$.
The homogeneous Laurent lift is defined by
\[
\operatorname{ch}(T)
=
\frac{P_B}{P_A}\widetilde{\ch}(T'),
\]
where $P_A$ and $P_B$ are the corresponding products
of trivial Pl\"ucker coordinates. The image of $\ch(T)$ under the specialization $f_a\mapsto1$ is
$\widetilde{\ch}(T')$. 
By \cite[Theorem~3.3]{DGL26},
$\ch(T)$ is in fact a polynomial in the Pl\"ucker coordinates, and
$\{\ch(T):T\in{\SSYT}(k,[n])\}$ is the dual canonical basis of
$\bC[\Gr_{k,n}]$.  We emphasize that the construction of
${\mathsf R}_{k,n}^{\bfy}$, $\bfy=(1,\dots,1)$ involves only the noncyclic frozen
variables $f_1,\ldots,f_{n-k+1}$, whereas the cyclic specialization
considered in this paper involves all $n$ frozen variables.

\subsection{Gelfand--Tsetlin degeneration}
\label{subsec:diagonal-sagbi-new}

We recall the Gelfand--Tsetlin toric degeneration of the Grassmannian.
Let $X=(x_{ij})_{i\in [k],j\in [n]}$ be a generic
$k\times n$ matrix and identify $R_{k,n}$ with the subalgebra of
$$\mathscr R:=\bC[x_{ij}:1\leq i\leq k,\ 1\leq j\leq n]$$ generated
by the maximal minors of $X$.  By a \emph{monomial order} we mean a
total well-order on monomials compatible with multiplication.  Fix the
lexicographic monomial order $\succ$ induced by
\[
x_{11}\succ x_{12}\succ\cdots\succ x_{1n}
\succ x_{21}\succ\cdots\succ x_{2n}
\succ\cdots\succ x_{k1}\succ\cdots\succ x_{kn}.
\]
For $0\neq h\in\mathscr R$, let
$\operatorname{LM}_{\prec}(h)$ denote its leading monomial.

If $A\subseteq\mathscr R$ is a subalgebra, its \emph{initial algebra}
is
$\operatorname{in}_{\prec}(A)
=\bC[\operatorname{LM}_{\prec}(h):0\neq h\in A]$.
A subset $\mathcal G\subseteq A$ is a \emph{SAGBI basis} if
$\operatorname{in}_{\prec}(A)$ is generated by the leading monomials
of the elements of $\mathcal G$.  Here SAGBI stands for
subalgebra analogue of Gr\"obner bases for ideals 
\cite{KaMa89,RS90}.  A finite SAGBI basis
admits a compatible weight filtration and hence a flat degeneration to
the initial algebra; see \cite[Chapter~11]{Stu96}.

For $I=(i_1<\cdots<i_k)$, define
$\ov P_I:=x_{1,i_1}x_{2,i_2}\cdots x_{k,i_k}$.  The determinant
expansion gives
$\mathrm{LM}_{\prec}(P_I)=\ov P_I.$
For a semistandard tableau $T=(I_1,\ldots,I_r)$, put
$\ov P_T:=\ov P_{I_1}\cdots\ov P_{I_r}$.  Equivalently,
\begin{equation}\label{eq:dT-row-multiplicities}
\ov P_T
=
\prod_{i=1}^k\prod_{p=1}^n
x_{i,p}^{\mu_{i,p}(T)},
\end{equation}
where $\mu_{s,p}(T)$ is equal to the number of times $p$ appears in the
$s$-th row of $T$.
Consequently, the map $T\mapsto\ov P_T$ is injective.  Since row
multiplicities are additive under row-wise union, one also has
$\ov P_{S\cup T}=\ov P_S\ov P_T.$

The Pl\"ucker coordinates are known to form a SAGBI basis for $R_{k,n}$ with
respect to $\prec$; see \cite{SZ93}.  In particular, 
$\{\ov P_T:T\in{\SSYT}(k,[n])\}$
is a $\bC$-basis of $\ov R_{k,n}$.
We note also that there is a nonnegative integral weight filtration
on $R_{k,n}$ such that the associated graded algebra is isomorphic to $\ov R_{k,n}$ and $P_I$ correspond to $\ov P_I$.
The corresponding Rees algebra $\mathcal R_{k,n}$ is flat over $\bC[t]$ and
satisfies
\[
\mathcal R_{k,n}/(t-1)\simeq R_{k,n},
\qquad
\mathcal R_{k,n}/(t)\simeq\ov R_{k,n}.
\]

The defining relations of  $R_{k,n}$ are the classical Pl\"ucker relations. 
The defining relations of $\ov R_{k,n}$ are also quadratic and are explicitly given by
\begin{equation}\label{eq:Hibi-relation}
\ov P_I\ov P_J
=
\ov P_{I\wedge J}\ov P_{I\vee J},
\end{equation}
where $\wedge$ and $\vee$ are standard meet and join operations in the distributive lattice
$L_{k,n}:=\binom{[n]}k$ with componentwise order
\[
(i_1<\cdots<i_k)\leq(j_1<\cdots<j_k)
\quad\Longleftrightarrow\quad
i_r\leq j_r\quad(1\leq r\leq k).
\]
Hence $\ov R_{k,n}$ is realized as the  Hibi algebra of $L_{k,n}$ \cite{Hibi87}.
Under the usual identification of $L_{k,n}$ with the order ideals of
the rectangle poset $[k]\times[n-k]$,
$\Spec(\ov R_{k,n})$ is the affine Hibi toric variety for the order polytope
\cite{Sta86}.

\subsection{Positroid varieties}
\label{subsec:positroid-varieties-prelim}

We recall the positroid geometry used later; see
\cite{Pos06,KLS13}.  For $V\in\Gr_{k,n}$, let
$\mathcal M(V):=\{I\in\tbinom{[n]}k:P_I(V)\neq0\}$, the
rank-$k$ matroid associated with $V$.  Let
$(\Gr_{k,n})_{\geq0}\subseteq\Gr_{k,n}(\bR)$ denote the totally
nonnegative Grassmannian, consisting of points represented by
$k\times n$ real matrices whose maximal minors are all
nonnegative.  A rank-$k$ matroid $\mathcal M$ on $[n]$ is
called a \emph{positroid} if
$\mathcal M=\mathcal M(V)$ for some
$V\in(\Gr_{k,n})_{\geq0}$. The positroids are in bijection with various combinatorial objects: juggling patterns, Grassmann necklaces, bounded affine permutations, etc. \cite{Pos06}. 

Now let us consider the complex picture.
For a positroid $\mathcal M$, we write
$\Pi^\circ_{\mathcal M}\subseteq\Gr_{k,n}$ for the corresponding
open positroid variety and
$\Pi_{\mathcal M}:=\overline{\Pi^\circ_{\mathcal M}}$ for its
closure.  The open positroid varieties form a stratification of
$\Gr_{k,n}$.

A positroid $\mathcal M$ is equivalently encoded by its Grassmann
necklace
$\mathcal I(\mathcal M)=(I_1(\mathcal M),\ldots,I_n(\mathcal M))$,
where $I_a(\mathcal M)$ is the lexicographically minimal basis of
$\mathcal M$ with respect to the cyclic order beginning at $a$.
For this cyclic order, $C_a$ is the smallest $k$-subset.  Hence
$I_a(\mathcal M)=C_a$ if and only if $C_a\in\mathcal M$.  Since
$f_a=P_{C_a}$, this is equivalent to $f_a$ being nowhere vanishing
on $\Pi^\circ_{\mathcal M}$.  Consequently, the simultaneous
nonvanishing of $f_1,\ldots,f_n$ gives the open positroid variety
whose Grassmann necklace is $(C_1,\ldots,C_n)$.

We shall use that $\Pi_{\mathcal M}$ is irreducible, normal, and
Cohen--Macaulay, and that its homogeneous ideal in the Pl\"ucker
coordinate ring is generated by the $P_I$ with
$I\notin\mathcal M$; see \cite{KLS13,Lam14}.

For the noncyclic family, we also use interval positroids.  A
\emph{rank set} in $[n]$ is a collection
$\mathcal R=\{[a_1,b_1],\ldots,[a_{n-k},b_{n-k}]\}$ of ordinary
intervals such that the $a_i$'s are pairwise distinct and the
$b_i$'s are pairwise distinct.  For an interval
$[i,j]\subseteq[n]$, set
\[
\rho_{\mathcal R}(i,j)
:=
j-i+1-
\#\{[a,b]\in\mathcal R:[a,b]\subseteq[i,j]\}.
\]
For $V\in\Gr_{k,n}$, let
$r_{[i,j]}(V):=\dim\pi_{[i,j]}(V)$, where
$\pi_{[i,j]}:\bC^n\to\bC^{j-i+1}$ is the coordinate projection.
The rank conditions
$r_{[i,j]}(V)\leq\rho_{\mathcal R}(i,j)$ for all
$1\leq i\leq j\leq n$ define a positroid
$\mathcal M_{\mathcal R}$, called an \emph{interval positroid}.
Its positroid variety
\[
\Pi_{\mathcal M_{\mathcal R}}
=
\{V\in\Gr_{k,n}:
r_{[i,j]}(V)\leq\rho_{\mathcal R}(i,j)
\text{ for all }1\leq i\leq j\leq n\}
\]
is called an \emph{interval positroid variety}; see
\cite{Kn14,KL14,Paw18}.  Thus interval positroids are precisely the
positroids whose defining rank conditions may be expressed using
ordinary nonwrapping intervals.  This distinction mirrors the two
families above: the cyclic zero fiber is governed by general
positroid varieties, whereas the noncyclic zero fiber is governed by
interval positroid varieties.

\section{Algebras: cyclic case}\label{sec:algcyc}

In this section, we study the algebras $R_{k,n}^{\bfx}$ and
$\ov R_{k,n}^{\bfx}$ obtained by specializing all $n$ frozen
Pl\"ucker coordinates to arbitrary scalars.
We construct bases for both families using noncrossing monomials
and cyclically reduced semistandard tableaux, and an additional
basis for $R_{k,n}^{\bfx}$ using cyclically reduced dual canonical
basis elements.

\subsection{Noncrossing world}
We first recall the notion of noncrossing tableaux from \cite{Pyl09,PPS10,SSW17}. 

Let us start with two collections $I,J\in \binom{[n]}{k}$,
$I=(i_1<\dots <i_k)$, $J=(j_1<\dots <j_k)$.
First, remove all common elements of $I$ and $J$,
namely the elements of $I\cap J$.
Let us denote the remaining parts by $I'=(i'_1<\dots <i'_s)$ and $J'=(j'_1<\dots <j'_s)$.
The collections $I$ and $J$ are said to be noncrossing if for all $a=1,\dots,s-1$,  
neither of the following two conditions holds: $i'_a<j'_a<i'_{a+1}<j'_{a+1}$,
$j'_a<i'_a<j'_{a+1}<i'_{a+1}$.

A rectangular tableau with $k$ rows is called noncrossing
if its columns are pairwise noncrossing.
A Pl\"ucker monomial $P_{I^1}\cdots P_{I^r}$
is called noncrossing if the columns
$I^1,\ldots,I^r$ are pairwise noncrossing.

One has the following Lemma.
\begin{lem}\label{lem:NC}
The noncrossing monomials form a basis of $R_{k,n}$. For any $a\in [n]$, the collection $C_a$
is noncrossing with any other collection $I$.
\end{lem}

The $n$ elements $C_a$ are called solid in \cite{PPS10} (see also \cite{SSW17}).

For $a\in [n]$, we define the elements $C_a\in\binom{[n]}{k}$ as follows:
\[
C_a = \begin{cases} (a,a+1,\dots,a+k-1), & a \le n-k+1,\\
(1,\dots,k-(n-a+1),a,\dots,n), & a > n-k+1
\end{cases}.
\]
So, in the second case, $a>n-k+1$, we start at $a$, move towards $n$, and
then restart from $1$ until we get $k$ numbers.

Let $f_a=P_{C_a}$ be the frozen variables in the Grassmannian cluster algebra
and let $S$ be the polynomial ring in variables $f_a$, $a\in [n]$.

For a collection $\bfx\in\bA^n$, let 
\[
R^\bfx_{k,n} = R_{k,n}/(f_1-x_1,\dots,f_n-x_n). 
\]
In particular, the zero fiber $R^0_{k,n}$ inherits a natural grading from $R_{k,n}$
counting the number of Pl\"ucker variables in a monomial. We write 
$R^0_{k,n}=\bigoplus_{m\ge 0} R^0_{k,n}(m)$.

\begin{lem}\label{lem:cycfree}
The algebra $S$ acts freely on $R_{k,n}$.    
\end{lem}
\begin{proof}
By Lemma~\ref{lem:NC}, the noncrossing monomials form a
$\bC$-basis of $R_{k,n}$, and every frozen column $C_a$
is noncrossing with every column. Every noncrossing monomial
therefore factors uniquely as a monomial in
$f_1,\ldots,f_n$ times a noncrossing monomial containing
no frozen factors. This gives an $S$-basis of $R_{k,n}$.
\end{proof}

\begin{cor}
For any $\bfx\in\bA^n$, the noncrossing monomials not divisible by any variable of the 
form $P_{C_a}$, $a\in [n]$, form a basis of $R^\bfx_{k,n}$.
\end{cor}

Let $H_{k,n}(t)$ be the Hilbert series of $R_{k,n}$. Hence, 
\[
H_{k,n}(t) = \sum_{m\ge 0} t^m \dim R_{k,n}(m) = \sum_{m\ge 0} t^m\dim V_{m\om_k}.
\]
We denote by $H^0_{k,n}=\sum_{m\ge 0} t^m \dim R^0_{k,n}(m)$ the Hilbert series 
of the algebra $R^0_{k,n}$.

\begin{cor}\label{cor:H0}
One has $H^0_{k,n}(t)=(1-t)^nH_{k,n}(t)$.
\end{cor}
\begin{proof}
Since $R_{k,n}$ is a free graded module over
$S=\bC[f_1,\ldots,f_n]$, with each $f_a$ of degree one,
we have $H_{k,n}(t)=\frac{H^0_{k,n}(t)}{(1-t)^n}$.
\end{proof}

\begin{example}
Let $n=4$, $k=2$. Then $R_{k,n}$ is a quotient of the polynomial algebra in six Pl\"ucker 
variables $P_{ij}$, $1\le i<j\le 4$, by the ideal generated by $P_{12}P_{34}-P_{13}P_{24}+P_{14}P_{23}$. 
Hence $H_{2,4}(t) = \frac{1-t^2}{(1-t)^6}$.

The zero-fiber algebra $R^0_{2,4}$ is  a quotient of the polynomial ring in variables $P_{13}$,
$P_{24}$ (other four variables are of the form $P_{C_a}$) modulo the relation $P_{13}P_{24}=0$.
The only noncrossing monomials in variables $P_{ij}$ that do not contain $P_{12}$, $P_{23}$,
$P_{34}$, $P_{14}$ are $P_{13}^m$, $P_{24}^m$ for $m>0$, and the empty one.
The Hilbert series is given by $H^0_{2,4}(t) = 1 + 2t + 2t^2 + \dots = \frac{1+t}{1-t}$.
\end{example}

\subsection{Noncrossing: the GT case}
We recall the following realization of the GT algebra $\ov R_{k,n}$. 

Let
$\eO_{k,n}$ be the order polytope \cite{Sta86} of the poset consisting of pairs $(i,j)$,
$1\le j\le k < i \le n$ with the poset structure defined by $(i,j)\le (i',j')$
if and only if $i\le i'$, $j\le j'$. In other words, $\eO_{k,n}$ is a polytope of 
$\bR^{k(n-k)}$ consisting of points $\bfs=(s_{i,j})_{i,j}$ such that 
$0\le s_{i,j}\le 1$ and $s_{ij}\le s_{i'j'}$ whenever $(ij)\le (i'j')$.

For $m\in\bZ_{\ge0}$, let $m\eO_{k,n}$
be the $m$-fold Minkowski sum of
$\eO_{k,n}$. The degree-$m$ component
$\ov R_{k,n}(m)$ is identified with
the span of the elements
\[
e^{(m,\bfs)},\qquad
\bfs\in m\eO_{k,n}\cap\bZ^{k(n-k)}.
\]
Multiplication is given by
\[
e^{(m,\bfs)}e^{(m',\bfs')}
=e^{(m+m',\bfs+\bfs')}.
\]
Since the order polytope has the integer
decomposition property, $\ov R_{k,n}$
is generated by the degree-one elements
$e^{(1,\bfs)}$, where
$\bfs\in\eO_{k,n}\cap\bZ^{k(n-k)}$.
There is a standard bijection between
these generators and the subsets
$I\in\binom{[n]}{k}$.
Hence $\ov R_{k,n}$ can be identified
with the quotient of the polynomial ring
in variables $\ov P_I$ by the Hibi ideal
described in
Section~\ref{subsec:diagonal-sagbi-new}.

Let $\ov S$ be the polynomial algebra generated by variables $\ov f_a = \ov P_{C_a}$.
We have the following proposition.

\begin{prop}
The algebra $\ov S$ acts freely on $\ov R_{k,n}$.     
For any $\bfx\in\bA^n$, the noncrossing monomials not divisible by any variable of the 
form $\ov P_{C_a}$, $a\in [n]$, form a basis of $\ov R^\bfx_{k,n}$.
\end{prop}
 
\begin{proof}
In \cite{PPS10} the authors prove that for any $\bfs\in m\eO_{k,n}\cap \bZ^{k(n-k)}$
the element $e^{(m,\bfs)}\in \ov R_{k,n}$ admits a unique decomposition
as a product of degree one elements $e^{(1,\bfs')}$ corresponding to a noncrossing
monomial in coordinates $\ov P_I$. This proves the second claim. For the first
claim, we again note that a product of $\ov P_{C_a}$ and a noncrossing monomial
is again noncrossing.
\end{proof} 

\subsection{Cyclically reduced tableaux}
Let us introduce a notion of cyclically reduced semistandard tableaux. 

Given a semistandard tableau $T$ of height $k$ with $m$ columns, we call
it cyclically reduced if there is no $a\in [n]$ such that 
$T=T'\cup C_a$
for a semistandard tableau $T'$ of width $m-1$, where $\cup$ is the row-wise union.

\begin{rem}
It is important to write the entries of $C_a$ in the increasing order, i.e.
for $a>n-k+1$ we write $C_a = (1,\dots, k-(n-a+1),a,\dots,n)$. 
\end{rem}

We call a monomial $P_T = \prod_{s=1}^m P_{I^{(s)}}$ cyclically reduced if the 
tableaux $T$ with rows $I^{(s)}$ is semistandard and cyclically reduced. 
The same definition is given for  $\ov P_T = \prod_{s=1}^m \ov P_{I^{(s)}}$.

\begin{lem}
Every semistandard tableau has a
unique decomposition
\[
T=T'\cup C_1^{m_1}\cup\cdots\cup C_n^{m_n},
\]
where $T'$ is cyclically reduced.
\end{lem}
\begin{proof}
By \cite{EL25} every $T\in \SSYT(k,[n])$ admits a unique unordered  decomposition into a union of 
pairwise noncrossing one-column tableaux. Since $C_a$ is noncrossing with any other column, we obtain the desired claim.
\end{proof}

For a collection $\bfx\in\bA^n$, we define 
\[
\ov R_{k,n}^{\bfx}
=
\ov R_{k,n}/(\ov f_1-x_1,\ldots,\ov f_n-x_n).
\]

\begin{thm}
Cyclically reduced monomials in variables $P_I$ (resp., $\ov P_I$) form a basis
of $R^\bfx_{k,n}$ (resp., $\ov R^\bfx_{k,n}$) for any $\bfx\in\bA^n$.
\end{thm}
\begin{proof}
Let us first consider the GT degenerate algebra $\ov R^\bfx_{k,n}$. Recall that 
any $T\in \SSYT(k,[n])$ admits a unique decomposition 
\begin{equation}\label{eq:decompcyc}
T = T'\cup C_1^{m_1}\cup\cdots\cup C_n^{m_n}, \qquad T'\in\SSYT^{cr}(k,[n]).
\end{equation}
Since $\ov P_T$ respects the row-wise union, one gets 
$\ov P_T = \ov P_{T'}\ov f_1^{m_1}\dots \ov f_n^{m_n}$. Hence cyclically reduced 
monomials in variables $\ov P_I$ form an $\ov S$ basis of $\ov R_{k,n}$.
They are also $\ov S$ linearly independent, because the decomposition \eqref{eq:decompcyc}
is unique and $\ov P_T$, $T\in\SSYT(k,[n])$ form a basis of $\ov R_{k,n}$. We conclude 
that $\ov R_{k,n} = \bigoplus_{T\in\SSYT^{cr}_{k,n}} \ov S\cdot \ov P_T$, which proves
the desired claim.

Now let us consider the action of the polynomial algebra $S$ on $R_{k,n}$. We consider 
the monomials $P_T$, $T\in\SSYT^{cr}_{k,n}$. Since $\ov P_T$ is the leading term
of $P_T$ and $\ov f_a$ is the leading term of $f_a$, we conclude that $P_T$, 
$T\in\SSYT^{cr}_{k,n}$ span $R_{k,n}$ as $S$ module. Now $S$ acts freely on $R_{k,n}$,
graded dimensions of $R_{k,n}$ and $\ov R_{k,n}$ coincide, and  $\ov P_T$, 
$T\in\SSYT^{cr}_{k,n}$ form an $\ov S$ basis of $\ov R_{k,n}$. We conclude that 
$R_{k,n} = \bigoplus_{T\in\SSYT^{cr}_{k,n}} S\cdot P_T$, which gives the 
desired claim.
\end{proof}

\begin{cor}
For any $m\geq 0$, the number of noncrossing tableaux
with $k$ rows and $m$ columns containing no cyclic
frozen columns equals the number of cyclically reduced
semistandard tableaux of the same shape.
\end{cor}
\begin{proof}
The two families give homogeneous bases of
$R^0_{k,n}$, so their cardinalities agree in each degree.
\end{proof}

\subsection{Dual canonical basis and cyclic frozen multiplication}

We first recall the dual canonical basis input used in this section.
By \cite{Ka91,Lus93} the Grassmannian
coordinate ring has a dual canonical basis $\{G_T:T\in \SSYT(k,[n])\}$.
We use the tableau parametrization for which the transition from the
standard monomial basis is unitriangular.  Recall (see 
section \ref{subsec:quantum-affine-algebras}) the homogeneous polynomials
$\ch(T)\in R_{k,n}$ satisfying $\ch(T)=G_T$.
Thus $\{\ch(T):T\in \SSYT(k,[n])\}$ is the dual canonical basis of $R_{k,n}$.

We need the following multiplicativity statement.

\begin{lem}\label{lem:GP}
Let $T\in \SSYT(k,[n])$, $F=C_a$ for some $a\in [n]$ 
written increasingly as a one-column tableau. Then
\[
G_TP_F=G_{T\cup F}.
\]
\end{lem}
\begin{proof}
Let $G_{q,T}$ be the quantum dual canonical basis element specializing to $G_T$, 
and let $\Delta_F$ be the quantum Pl\"ucker coordinate indexed by $F$.  A cyclic 
interval $F$ is weakly separated from every $k$-subset $I$.  By the
quasi-commutation theorem of Leclerc--Zelevinsky
\cite{LZ98},
\[
\Delta_F\Delta_I=q^{c(F,I)}\Delta_I\Delta_F,
\]
where $c(F,I)$ is an integer determined by $F,I$. For fixed $F$, the exponent $c(F,I)$ is additive in the Pl\"ucker factors and depends only on the Pl\"ucker degree and total content.
Since $G_{q,T}$ is homogeneous for these gradings, there is
$r\in\bZ$ such that
\[
\Delta_FG_{q,T}=q^rG_{q,T}\Delta_F.
\]

The element $\Delta_F$ is a quantum flag minor.  Caldero's
multiplicativity theorem for quantum flag minors
\cite{Cal03} therefore implies that, up to a power of $q$,
$G_{q,T}\Delta_F$ is a single quantum dual canonical basis element.
After specializing $q=1$, one gets $G_TP_F=G_U$ for some semistandard tableau $U$.

It remains to identify $U$.  In the usual tableau order for standard
monomial theory, the dual canonical basis is unitriangular over the
standard monomial basis, with highest term $P_T$, and this order is
compatible with adjoining a fixed column.  The standard straightening
law therefore shows that the unique highest standard monomial in
$G_TP_F$ is $P_{T\cup F}$, with coefficient $1$.  Hence
$U=T\cup F$, and
\[
G_TP_F=G_{T\cup F}.
\]
The statement for a trivial tableau follows by repeated multiplication
by the ordinary consecutive columns $C_a$.
\end{proof}

Recall that for $a\in\bZ_n$,  $C_a$ 
denotes the one-column semistandard tableau obtained by writing the
cyclic interval $\{a,a+1,\ldots,a+k-1\}\subset\bZ_n$
in increasing order as a $k$-subset of $[n]$.   A
semistandard tableau is called cyclically reduced if no $C_a$
divides it in the tableau monoid. We denote by $\SSYT^{cr}(k,[n])$
the set of cyclically reduced tableaux.  We use the unique cyclic
reduced decomposition:
\[
T = T'\cup C_1^{m_1}\cup\cdots\cup C_n^{m_n}, \qquad
T'\in\SSYT^{cr}(k,[n]).
\]

By Lemma~\ref{lem:GP} and the identification $G_T=\ch(T)$, one has 
\[
\ch(T)P_{C_a}=\ch(T\cup C_a).
\]
Recall the notation 
$S=\bC[f_1,\dots,f_n]$, $\ov S = \bC[\ov f_1,\dots,\ov f_n]$.

\begin{thm}
\label{thm:cyclic-ch-and-diagonal-bases}
The elements $\ch(T)$, $T\in\SSYT^{cr}(k,[n])$,
form an $S$-basis of $R_{k,n}$.
Moreover,
$\operatorname{LM}_{\prec}(\ch(T))=\ov P_T$.
These leading monomials form an $\ov S$-basis of
$\ov R_{k,n}$.
\end{thm}

\begin{proof}
By the unique cyclic reduced decomposition, every tableau has a unique
expression
\[
T=T'\cup C_1^{m_1}\cup\cdots\cup C_n^{m_n}
\]
with $T'$ cyclically reduced.  Repeated application of
Lemma~\ref{lem:GP} gives
\[
\ch(T)
=
\ch(T')f_1^{m_1}\cdots f_n^{m_n}.
\]
Since the $\ch(T)$'s form a basis of $R_{k,n}$, this proves the
first direct-sum decomposition.  

We are left to prove that the initial part of $\ch T$ is equal to $\ov P_T$.
Recall the embedding of the algebra $R_{k,n}$ into the polynomial
ring in variables $x_{i,j}$, $i\in [k]$, $j\in [n]$. 
The leading monomial of $P_I$ 
with respect to the lex ordering
is equal to $\ov P_I = \prod_{s=1}^k x_{s,i_s}$ and the leading term of $P_T$ is equal to
\[
\ov P_T = \prod_{s=1}^k\prod_{p=1}^n x_{s,p}^{\mu_{s,p}(T)},
\]
where $\mu_{s,p}(T)$ is equal to the number of times $p$ appears in the
$s$-th row of $T$.
By \cite{Br06}, the transition from the dual
canonical basis to the standard monomial basis
is unitriangular:
\[
\ch(T)=P_T+\sum_{B<T}c_BP_B,
\]
where $B<T$ implies $B\neq T$ and
\[
\sum_{i=1}^r\sum_{p=1}^t\mu_{i,p}(B)
\le
\sum_{i=1}^r\sum_{p=1}^t\mu_{i,p}(T)
\]
for all $1\le r\le k$ and $1\le t\le n$.

Let $(r,s)$ be the first position, in
row-major order, at which $T$ and $B$
differ. 
The preceding inequalities imply
$T_{r,s}<B_{r,s}$.
Indeed, if $B_{r,s}<T_{r,s}$, then the first
$r$ rows of $B$ would contain more entries
at most $B_{r,s}$ than those of $T$,
contradicting the preceding inequalities.
Since the entries in each row are weakly
increasing, the lexicographic order on
diagonal monomials gives
$\ov P_T\succ\ov P_B$.

Since the diagonal monomials $\ov P_T$ are
pairwise distinct and the $\ch(T)$'s form
a basis of $R_{k,n}$, the monomials
$\ov P_T$ form a basis of $\ov R_{k,n}$.
Moreover, the unique cyclic reduced decomposition
gives
$\ov P_T
=
\ov P_{T'}\ov f_1^{m_1}\cdots\ov f_n^{m_n}$.
Hence
\[
\ov R_{k,n}
=
\bigoplus_{T'\in\SSYT^{cr}(k,[n])}
\ov S\,\ov P_{T'}.
\]
\end{proof}

We summarize the results of this section in the following Corollary.

\begin{cor} \label{cor:5 bases cylic case}
Algebras $R^\bfx_{k,n}$ admit the following bases for all $\bfx\in\bA^n$:
\begin{enumerate}
\item noncrossing monomials $P_T$ without factors $P_{C_a}$, $a\in\bZ_n$,
\item cyclically reduced monomials $P_T$,
\item elements $\ch(T)$ for $T$ cyclically reduced.
\end{enumerate}    
Algebras $\ov R^\bfx_{k,n}$ admit the following bases for all $\bfx\in\bA^n$:
\begin{enumerate}
\item noncrossing monomials $\ov P_T$ without factors $\ov P_{C_a}$, $a\in\bZ_n$,
\item cyclically reduced monomials $\ov P_T$.
\end{enumerate}    
\end{cor}

\begin{example}
Let $(k,n)=(3,6)$. The frozen variables are
\[
f_1=P_{123},\quad f_2=P_{234},\quad
f_3=P_{345},\quad f_4=P_{456},\quad
f_5=P_{156},\quad f_6=P_{126}.
\]
Fix $\bfx=(x_1,\ldots,x_6)\in\bA^6$, so that
$f_a=x_a$ in $R_{3,6}^{\bfx}$. Let $\mathcal N$ be the set of finite multisets
of pairwise noncrossing $3$-subsets of $[6]$
containing none of the frozen columns
\[
123,\quad234,\quad345,\quad456,\quad
156,\quad126.
\]
We include the empty multiset.
Let $\mathcal T=\SSYT^{cr}(3,[6])$
be the set of cyclically reduced semistandard
tableaux, including the empty tableau.
By Corollary \ref{cor:5 bases cylic case}, the following
are three bases of $R_{3,6}^{\bfx}$:
\[
\begin{aligned}
\mathcal B_{\mathrm{nc}}
=\left\{\prod_{I\in M}P_I:
M\in\mathcal N\right\}, \ 
\mathcal B_{\mathrm{ss}}
=\{P_T:T\in\mathcal T\}, \ 
\mathcal B_{\mathrm{ch}}
=\{\ch(T):T\in\mathcal T\}.
\end{aligned}
\]

For example, consider the cyclically reduced
semistandard tableau $T=135\mid246$,
where $135$ and $246$ are its columns.
Then $P_{135}P_{246}\in\mathcal B_{\mathrm{ss}}$,
but this monomial is not noncrossing.
In the noncrossing monomial basis of $R_{3,6}^{\bfx}$,
we have
\[
P_{135}P_{246}
=P_{124}P_{356}+P_{145}P_{236}
+ x_3x_6+x_2x_5-x_1x_4.
\]
The dual canonical basis element $\ch(T)$
has the following expansions in the standard monomial
basis and the noncrossing monomial basis of $R_{3,6}$,
respectively:
\[
\begin{aligned}
\ch(T)
&=P_{135}P_{246}-P_{134}P_{256}
-P_{125}P_{346}+P_{124}P_{356}
-2P_{123}P_{456}\\
&=P_{145}P_{236}-P_{123}P_{456}.
\end{aligned}
\]
After specialization, the second expansion gives $\ch(T)=P_{145}P_{236}-x_1x_4$.

We now consider the Gelfand--Tsetlin
degeneration $\ov R_{3,6}^{\bfx}$.
Its two bases are
\[
\begin{aligned}
\ov{\mathcal B}_{\mathrm{nc}}
=\left\{\prod_{I\in M}\ov P_I:
M\in\mathcal N\right\},\quad 
\ov{\mathcal B}_{\mathrm{ss}}
=\{\ov P_T:T\in\mathcal T\}.
\end{aligned}
\]
Recall that $\ov P_{ijk}=x_{1i}x_{2j}x_{3k}$
for $i<j<k$, where the $x_{ij}$ are
matrix coordinates, distinct from the
specialization parameters $x_a$.
For the tableau $T=135\mid246$, the element $\ov P_{135}\ov P_{246}$ belongs to
$\ov{\mathcal B}_{\mathrm{ss}}$. We have
\[
\begin{aligned}
\ov P_{135}\ov P_{246}
=(x_{11}x_{23}x_{35})
  (x_{12}x_{24}x_{36})=(x_{11}x_{24}x_{35})
  (x_{12}x_{23}x_{36})=\ov P_{145}\ov P_{236}.
\end{aligned}
\]
Thus the element $\ov P_{135}\ov P_{246}$ is represented by
$\ov P_{145}\ov P_{236}$ in
$\ov{\mathcal B}_{\mathrm{nc}}$.
\end{example}

\section{Algebras: noncyclic case}\label{sec:algnoncyc}
In this section we restrict ourselves to variables 
$f_a=P_{C_a}$ and $\ov f_a = \ov P_{C_a}$ for $1\le a\le n-k+1$. 
In other words, we consider the nonwrapping sequences $a,a+1,\dots,a+k-1$ such that $a+k-1\le n$.
We denote the polynomial algebra in variables $f_a$ (resp., $\ov f_a$), $a\in [n-k+1]$ by ${\mathsf S}$ (resp., $\ov {\mathsf S}$). 

For $\bfy=(y_1,\ldots,y_{n-k+1})\in\bA^{n-k+1}$, define
\[
{\mathsf R}^\bfy_{k,n}
=R_{k,n}/(f_1-y_1,\ldots,f_{n-k+1}-y_{n-k+1}),
\]
\[
\ov{\mathsf R}^\bfy_{k,n}
=\ov R_{k,n}/(\ov f_1-y_1,\ldots,\ov f_{n-k+1}-y_{n-k+1}).
\]

We construct bases of ${\mathsf R}^\bfy_{k,n}$ and
$\ov{\mathsf R}^\bfy_{k,n}$ using noncrossing monomials
without frozen factors and reduced semistandard tableaux.
For ${\mathsf R}^\bfy_{k,n}$, we also obtain a basis
consisting of reduced dual canonical basis elements.

The following proposition is proved along the same lines as Lemmas \ref{lem:NC} and \ref{lem:cycfree}. 

\begin{prop}
The algebra $R_{k,n}$ (resp., $\ov R_{k,n}$)
is a free module over ${\mathsf S}$
(resp., $\ov{\mathsf S}$), with a basis consisting
of noncrossing monomials $P_T$ (resp., $\ov P_T$)
without factors $P_{C_a}$ (resp., $\ov P_{C_a}$),
$a\in[n-k+1]$.
\end{prop}

A semistandard tableau is called  reduced if no $C_a$, $a\in [n-k+1]$
divides it in the tableau monoid. We denote by $\SSYT^{r}(k,[n])$
the set of reduced tableaux.

By \cite{CDFL20}, we know that $T=T_{\rm red} \cup T''$ where $T''$'s columns are trivial. The explicit expression of $T''$ is as follows. 

For $T\in \SSYT(k,[n])$, $a\in [n-k+1]$, and $s\in [k]$, denote by $m_{a,s}(T)$ the number of 
entries of $T$ in the $s$-th row that are equal to $a+s-1$ and let $m_a(T)$ be the minimum 
of all $m_{a,s}(T)$.

\begin{lem} \label{lem:factorization of semistandard tableau factor out noncyclic frozen}
Every semistandard tableau has a
unique reduced decomposition
\[
T=T'\cup C_1^{m_1}\cup\cdots\cup C_{n-k+1}^{m_{n-k+1}},
\]
with $T'$ reduced.  
\end{lem}
\begin{proof} 
One can decompose
\[
T = T' \cup C_1^{m_1(T)} \cup\dots \cup C_{n-k+1}^{m_{n-k+1}(T)}.
\]
It remains to prove that $T'\in \SSYT(k,[n])$. The latter is equivalent to the following set of
inequalities: for any $s<k$ and any $r<n$, the number of entries not exceeding $r$ in the 
$s$-th row is at least the number of entries not exceeding $r+1$ in the $(s+1)$-st row.
Since removing a factor  $C_a$, $a\in [n-k+1]$ does not change the difference of these two
numbers, the condition $T\in \SSYT(k,[n])$ implies $T'\in \SSYT(k,[n])$.

For each $a$, the minimum defining $m_a(T)$ is attained
in at least one row. After removing the indicated frozen
columns, that row contains no entry $a+s-1$ in row $s$;
thus $T'$ is reduced. Conversely, if
$T=U\cup\prod_{a=1}^{n-k+1}C_a^{r_a}$ with $U$ reduced,
then $m_a(T)=m_a(U)+r_a=r_a$ for every $a$.
This proves uniqueness.
\end{proof}

We have the following theorem.

\begin{thm} \label{thm:5 bases noncyclic case}
The algebra $R_{k,n}$ admits the following
${\mathsf S}$-bases:
\begin{enumerate}
\item noncrossing monomials $P_T$ without factors $P_{C_a}$, $a\in [n-k+1]$,
\item reduced monomials $P_T$,
\item elements $\ch(T)$ for $T$ reduced.
\end{enumerate}    
The algebra $\ov R_{k,n}$ admits the following
$\ov{\mathsf S}$-bases:
\begin{enumerate}
\item noncrossing monomials $\ov P_T$ without factors $\ov P_{C_a}$, $a\in [n-k+1]$,
\item reduced monomials $\ov P_T$.
\end{enumerate}    
\end{thm}

\begin{proof}
The proof is similar to that of the cyclic case in
Section~\ref{sec:algcyc}, with the difference that we only
consider the nonwrapping columns $C_a$, $a\in[n-k+1]$.

First, each $C_a$ is noncrossing with every other column.
Hence every noncrossing monomial admits a unique factorization
into a noncrossing monomial without frozen factors and a
monomial in $f_1,\ldots,f_{n-k+1}$.
Since noncrossing monomials form bases of $R_{k,n}$ and
$\ov R_{k,n}$, this proves the first assertions.

Next, by Lemma \ref{lem:factorization of semistandard tableau factor out noncyclic frozen}, every semistandard tableau
has a unique decomposition
\[
T=T'\cup C_1^{m_1}\cup\cdots\cup
C_{n-k+1}^{m_{n-k+1}},
\]
where $T'$ is reduced. Since multiplication in $\ov R_{k,n}$ respects row-wise
union, we have
\[
\ov P_T=\ov P_{T'}\ov f_1^{m_1}\cdots
\ov f_{n-k+1}^{m_{n-k+1}}.
\]
The uniqueness of the decomposition therefore gives
\[
\ov R_{k,n}
=\bigoplus_{T'\in\SSYT^r(k,[n])}
\ov{\mathsf S}\cdot\ov P_{T'}.
\]
For $R_{k,n}$, recall that the leading monomial of $P_T$
is $\ov P_T$ and that of $f_a$ is $\ov f_a$.
Thus the elements
\[
P_{T'}f_1^{m_1}\cdots f_{n-k+1}^{m_{n-k+1}},
\qquad T'\in\SSYT^r(k,[n]),
\]
have distinct leading monomials $\ov P_T$, which form
a basis of the initial algebra $\ov R_{k,n}$.
The same leading-term argument as in
Section~\ref{sec:algcyc} shows that these elements
form a basis of $R_{k,n}$. Hence
\[
R_{k,n}
=\bigoplus_{T'\in\SSYT^r(k,[n])}
{\mathsf S}\cdot P_{T'}.
\]

Finally, by Lemma~\ref{lem:GP}, multiplication by
a frozen Pl\"ucker coordinate satisfies
\[
\ch(T)f_a=\ch(T\cup C_a).
\]
Consequently, the unique reduced decomposition gives
\[
\ch(T)=\ch(T')f_1^{m_1}\cdots
f_{n-k+1}^{m_{n-k+1}}.
\]
Since the elements $\ch(T)$ form the dual canonical
basis of $R_{k,n}$, we obtain
\[
R_{k,n}
=\bigoplus_{T'\in\SSYT^r(k,[n])}
{\mathsf S}\cdot\ch(T').
\]
This completes the proof.
\end{proof}

\begin{example}
Let $(k,n)=(3,6)$. The nonwrapping frozen
variables are
\[
f_1=P_{123},\quad f_2=P_{234},\quad
f_3=P_{345},\quad f_4=P_{456}.
\]
Fix $\bfy=(y_1,\ldots,y_4)\in\bA^4$,
so that $f_a=y_a$ in
${\mathsf R}_{3,6}^{\bfy}$.
Notice that $P_{156}$ and $P_{126}$
are not specialized in this family. Let $\mathcal N$ be the set of finite
multisets of pairwise noncrossing
$3$-subsets of $[6]$ containing none
of the columns
\[
123,\quad234,\quad345,\quad456.
\]
We include the empty multiset.
Let $\mathcal T=\SSYT^r(3,[6])$ be
the set of reduced semistandard tableaux,
including the empty tableau.
By Theorem \ref{thm:5 bases noncyclic case}, the following
are three ${\mathsf S}$-bases of $R_{3,6}$:
\[
\begin{aligned}
\mathcal B_{\mathrm{nc}}
=\left\{\prod_{I\in M}P_I:
M\in\mathcal N\right\},\ 
\mathcal B_{\mathrm{ss}}
=\{P_T:T\in\mathcal T\},\
\mathcal B_{\mathrm{ch}}
=\{\ch(T):T\in\mathcal T\}.
\end{aligned}
\]
Their images give bases of
${\mathsf R}_{3,6}^{\bfy}$ for all $\bfy$.
Consider the reduced semistandard tableau
$T=135\mid246$, where $135$ and $246$
are its columns.
Then $P_{135}P_{246}\in\mathcal B_{\mathrm{ss}}$,
but this monomial is not noncrossing.
By the Pl\"ucker relations, its expansion
in the noncrossing monomial basis of
${\mathsf R}_{3,6}^{\bfy}$ is
\[
\begin{aligned}
P_{135}P_{246}
=P_{124}P_{356}+P_{145}P_{236}+y_3P_{126}+y_2P_{156}-y_1y_4.
\end{aligned}
\]
Here $P_{124}P_{356}$,
$P_{145}P_{236}$, $P_{126}$,
$P_{156}$, and $1$ belong to
$\mathcal B_{\mathrm{nc}}$.
The corresponding dual canonical basis
element has the noncrossing expansion
$\ch(T)=P_{145}P_{236}-P_{123}P_{456}$
in $R_{3,6}$. Consequently, in
${\mathsf R}_{3,6}^{\bfy}$, we have
\begin{align*}
\ch(T)=P_{145}P_{236}-y_1y_4.
\end{align*}

We now consider the Gelfand--Tsetlin
degeneration
$\ov{\mathsf R}_{3,6}^{\bfy}$.
Its two bases are
\[
\begin{aligned}
\ov{\mathcal B}_{\mathrm{nc}}
=\left\{\prod_{I\in M}\ov P_I:
M\in\mathcal N\right\}, \quad 
\ov{\mathcal B}_{\mathrm{ss}}=\{\ov P_T:T\in\mathcal T\}.
\end{aligned}
\]
Recall that $\ov P_{ijk}=x_{1i}x_{2j}x_{3k}$
for $i<j<k$, where the $x_{ij}$ are
the matrix coordinates. For $T=135\mid246$, the element
$\ov P_{135}\ov P_{246}$ belongs to
$\ov{\mathcal B}_{\mathrm{ss}}$. We have
\[
\begin{aligned}
\ov P_{135}\ov P_{246}
=(x_{11}x_{23}x_{35})
  (x_{12}x_{24}x_{36})
=(x_{11}x_{24}x_{35})
  (x_{12}x_{23}x_{36})
=\ov P_{145}\ov P_{236}.
\end{aligned}
\]
Thus the element $\ov P_{135}\ov P_{246}$
is represented by $\ov P_{145}\ov P_{236}$
in $\ov{\mathcal B}_{\mathrm{nc}}$.
\end{example}

\section{Varieties: cyclic case}\label{sec:varcyc}

We study the geometry of the family $\Spec(R_{k,n}^{\bfx})$
and describe the irreducible components of the fiber
over $\bfx=0$ in terms of positroid varieties. We use standard algebro-geometric notions and techniques, see e.g. \cite{V25}.

Let $r = k(n-k)+1-n$. Let $Z_{k,n}^\bfx$ be the spectrum of the algebra $R^\bfx_{k,n}$. So, 
$Z_{k,n}^\bfx$ is an affine algebraic variety in $\binom{n}{k}$-dimensional affine space. 
We start with the following theorem describing basic properties of all fibers.
\begin{thm}
The family $Z_{k,n}^\bfx$ is flat over $\bA^n$. 
For any $\bfx\in\bA^n$, $Z_{k,n}^\bfx$ is nonempty, connected, and reduced, 
and has pure dimension $r$.
\end{thm}

\begin{proof}
By Lemma~\ref{lem:cycfree}, $R_{k,n}$ is a free
module over $S=\bC[f_1,\ldots,f_n]$.
Hence the family is flat, and every fiber
is nonempty.

Let $J_{k,n}$ be the Pl\"ucker ideal.
By \cite{SSW17}, there is a term order for
which $J_{k,n}$ has a Gr\"obner basis
whose leading monomials are the products
$P_IP_J$ with $I,J$ crossing.
No frozen variable $P_{C_a}$ divides
any of these leading monomials.

For $\bfx\in\bA^n$, set
\[
J_{k,n}^{\bfx}
=J_{k,n}+(P_{C_a}-x_a:a\in[n]).
\]
The leading monomial of $P_{C_a}-x_a$
is $P_{C_a}$, which is relatively prime
to every leading monomial of the
chosen Gr\"obner basis of $J_{k,n}$.
By Buchberger's product criterion
\cite[Chapter~2, Section~9]{CLO15},
adjoining the polynomials $P_{C_a}-x_a$
gives a Gr\"obner basis of $J_{k,n}^{\bfx}$.
Consequently,
\[
\operatorname{in}(J_{k,n}^{\bfx})
=
\operatorname{in}(J_{k,n})
+(P_{C_1},\ldots,P_{C_n}).
\]
This is a squarefree monomial ideal.
Therefore $J_{k,n}^{\bfx}$ is radical,
and every fiber is reduced.

To prove connectedness, equip
$R_{k,n}^{\bfx}$ with the filtration
induced by the Pl\"ucker degree.
The homogeneous free $S$-basis
constructed in Lemma~\ref{lem:cycfree}
gives an isomorphism
$\operatorname{gr}R_{k,n}^{\bfx}
\simeq R_{k,n}^{0}$.
The latter algebra is reduced and has
degree-zero part equal to $\bC$.
If an idempotent of $R_{k,n}^{\bfx}$
had positive filtration degree, its
nonzero leading term would be
nilpotent, a contradiction.
Thus every idempotent belongs to
$\bC$ and is either $0$ or $1$.
Hence every fiber is connected.

Finally, since $R_{k,n}$ is free over $S$,
the elements $f_a-x_a$, $a\in[n]$,
form a regular sequence.
As $R_{k,n}$ is Cohen--Macaulay
of dimension $k(n-k)+1$, every fiber
is of pure dimension $k(n-k)+1-n=r$.
\end{proof}

\subsection{Generic fiber}\label{sec:cycgenfib}
In this subsection we describe  basic properties of a generic fiber. 
\begin{thm}\label{thm:cycgen}
There exists an open subset $U\subset (\bA\setminus \{0\})^n$ such that for
any $\bfx\in U$ the fiber $Z_{k,n}^\bfx$ is smooth and irreducible. It contains 
an open subset isomorphic to $(\bC^\times)^r$ and thus is rational.
\end{thm}
\begin{proof}
Let $\widehat{\Gr}_{k,n} = \mathrm{Spec}(R_{k,n})$ be the cone over the Grassmannian $\Gr_{k,n}$.
By \cite{Sc06}, $\mathrm{Spec}(R_{k,n})$ contains an open chart  
$U\simeq (\bC^\times)^r\times (\bC^\times)^n$. More precisely, one chooses 
a cluster with mutable variables $z_1,\dots,z_r$. Then $U$ is equal to the set
of points where all variables $z_i$ and all frozen variables $f_a$ do not vanish.
The identification $U\simeq (\bC^\times)^r\times (\bC^\times)^n$ is given by 
evaluating $z_i$ and $f_a$. Hence for each point $\bfx$, $x_i\ne 0$, one gets 
an open subset $U\cap Z_{k,n}^\bfx$ of the fiber.

Now let us consider $\widehat{\Gr}^\circ_{k,n}= \widehat{\Gr}_{k,n}\setminus \{0\}$. 
It is smooth as the
complement of the zero section in the  tautological line bundle over the smooth
Grassmannian.
Let $\pi: \widehat{\Gr}_{k,n}\to \bC^n$
be the projection map, $p\mapsto (f_1(p),\dots,f_n(p))$. Every fiber of $\pi$ over a point
of $(\bC^\times)^n$ sits inside $\widehat{\Gr}_{k,n}^\circ$. 
The characteristic-zero generic smoothness theorem now implies that there is a nonempty 
Zariski-open $B\subset (\bC^\times)^n$  such that the whole map
$ \pi^{-1}(B)\to B$ is smooth. Hence every fiber over $B$ is smooth; surjectivity of 
$\pi$ ensures these fibers are nonempty. 

We are left to show that the torus $(\bC^\times)^r$ is dense in every fiber over the set $B$
(in particular, this would imply irreducibility).
Recall the torus $U\simeq (\bC^\times)^r\times (\bC^\times)^n\subset \widehat{\Gr}_{k,n}$.
Put $M = \widehat{\Gr}_{k,n}\setminus U$. Since $\widehat{\Gr}_{k,n}$ is irreducible and $U$ is nonempty and open, 
$\dim M\leq r+n-1$. By the generic fiber dimension theorem, after shrinking the parameter open 
set, every fiber $M\cap Z_{k,n}^\bfx$  has dimension at most $r-1$. But every component of 
a smooth fiber $Z_{k,n}^\bfx$ has dimension $r$. Thus no component can lie in $M$: 
$U\cap Z_{k,n}^\bfx$ is dense. Since that open subset is the irreducible torus 
$(\bC^\times)^r$, the whole fiber is irreducible.
\end{proof}

\subsection{The special fiber}
In this subsection we consider the special fiber, i.e. $\bfx=0$. In this case  $R^\bfx_{k,n}$ inherits the grading from the Pl\"ucker
algebra and hence one can consider its projective spectrum. We define
\[
Y_{k,n} = \Proj(R^0_{k,n}) \subset \bP^{\binom{n}{k}-1},\quad 
Z_{k,n}^{0} = \mathrm{Cone}(Y_{k,n}).
\]
By definition, 
$Y_{k,n}$ is a subvariety of $\Gr_{k,n}$ cut out by the vanishing of the
cyclic Pl\"ucker coordinates $f_a$, $a\in [n]$. 
From the general discussion above, we know that 
$Y_{k,n}$ is a reduced equidimensional complete intersection (in $\Gr_{k,n}$) 
of dimension $k(n-k)-n$.

\begin{lem}\label{lem:openpositr}
$Y_{k,n}$  is the union of all open positroid varieties $\Pi^\circ_{\mathcal M}$  
such that $C_a\notin\mathcal M$ for all $a\in [n]$. 
\end{lem}
\begin{proof}
Recall that $C_a\in\mathcal M$ if and only if $f_a$ does not vanish on the open
positroid variety $\Pi^\circ_{\mathcal M}$. This implies that 
$\Pi^\circ_{\mathcal M}\subset Y_{k,n}$ whenever $C_a\notin \mathcal M$ for all $a$.
To complete the proof, we recall that open positroid varieties form a stratification of
the Grassmannian. Since $\Pi^\circ_{\mathcal M}\cap Y_{k,n}=\emptyset$ if $C_a\in 
\mathcal M$ for some $a$, we arrive at the desired equality. 
\end{proof}

\begin{thm}
$Y_{k,n}$ is the union of all codimension $n$ positroid varieties $\Pi_{\mathcal M}$  
such that $C_a\notin\mathcal M$ for all $a\in [n]$. Each such positroid variety
is an irreducible component of $Y_{k,n}$.
\end{thm}
\begin{proof}
We know that $\dim Y_{k,n}=k(n-k)-n$ and $Y_{k,n}$ is closed and equidimensional. 
Now our theorem follows from Lemma \ref{lem:openpositr}.    
\end{proof}

\begin{cor}
Let us consider the subposet of positroids $\mathcal M$ such that 
$C_a\notin\mathcal M$ for all $a\in [n]$. Then all the maximal elements of this 
subposet are of codimension $n$.
\end{cor}

\begin{example}
Let $(k,n)=(2,4)$. The fiber
$Z_{k,n}^\bfx$, $\bfx=(x_1,x_2,x_3,x_4)$ is the affine spectrum of the algebra 
\[
\bC[P_{13},P_{24}]/(P_{13}P_{24} - x_1x_3 - x_2x_4).
\]
The open subset $U$ from Theorem \ref{thm:cycgen} thus consists of $\bfx$ such that $x_1x_3+x_2x_4\ne 0$.
The special fiber over $\bfx=0$ is a cone over the zero-dimensional variety $\Proj(\bC[P_{13},P_{24}]/(P_{13}P_{24}))$ consisting of two points. 
\end{example}

\begin{example}
Let $(k,n)=(2,5)$. The fiber
$Z_{k,n}^\bfx$, $\bfx=(x_1,x_2,x_3,x_4,x_5)$ is the affine spectrum of the quotient of the polynomial algebra in variables $P_{13},P_{14},P_{24},P_{25},P_{35}$ by the ideal generated by five Pl\"ucker relations  
\begin{align*}
P_{13}P_{24} &= x_2P_{14}+x_1x_3,\\
P_{13}P_{25} &= x_1P_{35}+x_2x_5,\\
P_{14}P_{25} &= x_5P_{24}+x_1x_4,\\
P_{14}P_{35} &= x_4P_{13}+x_3x_5,\\
P_{24}P_{35} &= x_3P_{25}+x_2x_4.
\end{align*}
The special fiber over $\bfx=0$ is a cone over a one-dimensional variety with five irreducible components $K_a$, $1\le a\le 5$.
Each component is a line $\bP^1$; consecutive lines in the cyclic order
meet at a coordinate point; all other pairs are disjoint. The five intersection points are ordinary nodes.
\end{example}

\section{Varieties: noncyclic case}\label{sec:varnoncyc}

We study the geometry of the family $\Spec({\mathsf R}_{k,n}^{\bfy})$
obtained by specializing the nonwrapping frozen Pl\"ucker coordinates.
We show that the fiber over $\bfy=0$ is an affine cone over a
projective variety whose irreducible components are interval
positroid varieties of codimension $n-k+1$ in $\Gr_{k,n}$.

Let $s=(k-1)(n-k)$. We note that $s=k(n-k)+1 - (n-k+1)$. For $\bfy\in \bA^{n-k+1}$, let ${\mathsf Z}_{k,n}^\bfy= \Spec({\mathsf R}_{k,n}^{\bfy})$. 
The same argument as in the cyclic case gives the following proposition.
\begin{prop}
The family ${\mathsf Z}_{k,n}^{\bfy}$ is flat over
$\bA^{n-k+1}$. For every $\bfy\in\bA^{n-k+1}$,
the fiber ${\mathsf Z}_{k,n}^{\bfy}$ is nonempty,
connected, reduced, and of pure dimension $s$.
\end{prop}

\subsection{Generic fibers}
The situation with the generic fibers is much simpler in the nonwrapping case compared
to the cyclic case from section \ref{sec:cycgenfib}. 
\begin{prop}
Let $y_i\ne 0$ for all $i\in [n-k+1]$. Then the fiber ${\mathsf Z}_{k,n}^\bfy$ is isomorphic 
to $\bA^s$. Fibers ${\mathsf Z}_{k,n}^\bfy$ with all $y_i\ne 0$ form a  trivial vector
bundle over $\bA^{n-k+1} \setminus \bigcup_{i=1}^{n-k+1}\{y_i=0\}$.
\end{prop}
\begin{proof}
Let us construct an isomorphism ${\mathsf Z}_{k,n}^\bfy\simeq \bA^s$. 
Let $V\in \widehat{\Gr}_{k,n}$ be a point in a cone over $\Gr_{k,n}$. Since $f_1(V)=y_1\ne 0$,
$V$ can be represented by a $k\times n$ matrix $(v_1,\dots,v_n)$ (each $v_i$ is a 
$k$-vector) such that $v_1, \dots,v_k$ are basis vectors and 
$\det(v_1,\dots,v_k)=y_1$. Now we want to determine other $k$-vectors $v_{k+1},\dots,v_n$.
Suppose we have already fixed all vectors $v_i$ with $i\le t$. Since 
$\det(v_{t-k+1},\dots,v_t)=y_{t-k+1}\ne 0$, these vectors form a basis of the 
$k$-dimensional space and hence $v_{t+1}$ can be expressed as linear combination
of $cv_{t-k+1}+\sum_{i=2}^k \alpha_{t,i} v_{t-k+i}$ for some complex coefficients 
$\alpha_{t,i}$. We know that 
\[
\det(v_{t-k+1},\dots,v_t)=y_{t-k+1},\ \det(v_{t-k+2},\dots,v_t,v_{t+1})=y_{t-k+2}.
\]
Hence $(-1)^{k-1}c=y_{t-k+2}/y_{t-k+1}$ and all the coefficients $\alpha_{t,i}$ are 
free parameters. Since $k\le t <n$ and $2\le i\le k$, we obtain the identification 
${\mathsf Z}_{k,n}^\bfy\simeq \bA^s$ (recall $s=(n-k)(k-1)$)
for all $\bfy\in\bA^{n-k+1}$ with all $y_i$ being nonzero.
\end{proof}

\begin{example}
Let $(k,n)=(2,4)$ and consider the noncyclic
family with frozen variables
\[
f_1=P_{12},\qquad f_2=P_{23},\qquad
f_3=P_{34}.
\]
Write $u=P_{13}$, $v=P_{14}$, and $w=P_{24}$.
The Pl\"ucker relation gives
\[
{\mathsf R}_{2,4}^{\bfy}
\simeq
\bC[u,v,w]/(uw-y_1y_3-y_2v).
\]
If $y_2\ne0$, then $v$ is uniquely determined
by $u,w$, and hence
${\mathsf Z}_{2,4}^{\bfy}\simeq\bA^2$.
For $\bfy=0$, we have
${\mathsf R}_{2,4}^{0}
\simeq\bC[u,v,w]/(uw)$.
Consequently, the projective special fiber is
${\mathsf Y}_{2,4}
=V(uw)\subset\bP^2$,
which is the union of two projective lines.
These are interval positroid varieties.
For comparison, the cyclic zero fiber has
coordinate ring
$R_{2,4}^{0}\simeq\bC[u,w]/(uw)$,
and its projectivization consists of two points.
\end{example}

\subsection{The special fiber}
Let us consider the case $\bfy=0$. As in the cyclic case, 
${\mathsf Z}_{k,n}^0$ is a cone over a projective variety, which we denote by 
${\mathsf Y}_{k,n}\subset \bP^{\binom{n}{k}-1}$. Explicitly, 
${\mathsf Y}_{k,n} = \Proj({\mathsf R}_{k,n}^{0})$.

By definition, ${\mathsf Y}_{k,n}$ is a subvariety of $\Gr_{k,n}$ cut out by 
the vanishing of the nonwrapping Pl\"ucker coordinates $f_a$, $a\in [n-k+1]$.
Explicitly, if $V\in \Gr_{k,n}$ is represented by a matrix with columns $(v_1,\ldots,v_n)$
(where $v_i$ are $k$-vectors), then $V\in {\mathsf Y}_{k,n}$ if and only if 
$$\det(v_a,\dots,v_{a+k-1}) = 0\quad \text{ for all }\quad  a\in [n-k+1].$$
We know that ${\mathsf Y}_{k,n}$ is reduced and equidimensional of dimension $s-1$. 

\begin{thm}
${\mathsf Y}_{k,n}$ is the union of all open positroid varieties $\Pi^\circ_{\mathcal M}$  
such that $C_a\notin\mathcal M$ for all $a\in [n-k+1]$. Irreducible components of
${\mathsf Y}_{k,n}$ are the closed positroid varieties corresponding to maximal positroids 
satisfying $C_a\notin\mathcal M$ for all $a\in [n-k+1]$.
\end{thm}
\begin{proof}
Represent a point  $V\in \Gr_{k,n}$ by a matrix with columns $(v_1,\ldots,v_n)$ and set
\[
r_{ij}(V)=\mathrm{rank}(v_i,\ldots,v_j), \qquad i\le j.
\]
The variety ${\mathsf Y}_{k,n}$  is defined by conditions $f_a=0$, $a\in [n-k+1]$, i.e.
\[
r_{a,a+k-1}(V)\le k-1, \qquad 1\le a\le n-k+1.
\]

For a realized collection of ranks ${\bf r}=(r_{ij})$, consider
\begin{gather*}
\Pi_{\bf r}^\circ=\{V: r_{ij}(V) = r_{ij}\text{ for every interval }[i,j]\},\\
\Pi_r=\{V: r_{ij}(V)\le r_{ij}\text{ for every }[i,j]\}.
\end{gather*}
It is proved in \cite{Kn14} that $\Pi_{\bf r}^\circ$ is irreducible and that its closure 
is exactly $\Pi_r$ (these are the interval positroid varieties). 

Now take an irreducible component $X$ of ${\mathsf Y}_{k,n}$. 
We define 
\[
r_{ij}=\max_{V\in X}\mathrm{rank}(v_i,\ldots,v_j).
\]
The maximum is attained on a nonempty open subset of $X$. There are only finitely 
many intervals $[i,j]$, and $X$ is irreducible, so on a dense open subset of $X$, 
all the ranks have maximal constant values ${\bf r}=(r_{ij})$. Therefore $X\subset \Pi_{\bf r}$.

On this subset the length-$k$ ranks $r_{a,a+k-1}$ are at most $(k-1)$; the defining inequalities 
for $\Pi_{\bf r}$ imply the same throughout $\Pi_{\bf r}$. Hence
\[
X\subseteq\Pi_{\bf r}\subseteq {\mathsf Y}_{k,n}.
\]
Since $\Pi_{\bf r}$ is irreducible and $X$ is an irreducible component of ${\mathsf Y}_{k,n}$, 
we get $X=\Pi_{\bf r}$.  
\end{proof}

We would like to give a more explicit description of the irreducible components
in terms of the interval positroid varieties. 
The interval positroid varieties can be labeled by rank sets \cite{Paw18,KL14}.
A rank set is a collection of segments $[a_i,b_i]\subset [1,n]$, $i=1,\dots,n-k$ 
such that $a_1<\dots <a_{n-k}$ and all $b_i$ are pairwise distinct.
The expression $\sum_{i=1}^{n-k} (b_i-a_i) -  |\{i<j:\ b_i>b_j\}|$ computes the dimension of the corresponding positroid variety.

Let $F^{(k)}_t$ be the number of compositions of $t$
with parts no larger than $k$.
\begin{rem}
These numbers satisfy the $k$-step Fibonacci recurrence
\[
F_t^{(k)}=\sum_{j=1}^{k}F_{t-j}^{(k)}
\qquad (t\ge 1),
\]
with initial conditions
$F_0^{(k)}=1$ and $F_t^{(k)}=0$ for $t<0$.
\end{rem}

\begin{prop}
${\mathsf Y}_{k,n}$ is equal to the union of all closed interval positroid varieties
of dimension $(k-1)(n-k)-1$ contained in the zero locus of all $f_a$, $a\in [n-k+1]$.
The number of these positroid varieties is equal to $F^{(k)}_{n-k+1}-1$.
\end{prop}

\begin{proof}
It is proved in \cite{HS04} that the number of irreducible components of the variety of $k\times n$ matrices
whose consecutive maximal minors vanish is $F^{(k)}_{n-k+1}$. Exactly one
component consists of matrices of rank less than $k$, so it contributes no point to the Grassmannian.
Since the full-rank locus maps to ${\mathsf Y}_{k,n}$ as a principal $\mathrm{GL}_k$-bundle with irreducible fibers,
each remaining component gives a distinct irreducible component of ${\mathsf Y}_{k,n}$.

The paper \cite{HS04} gives the explicit description of the irreducible components. Namely, 
we take a sequence $\Gamma=([a_1,b_1],\ldots,[a_s,b_s])$
such that 
\[
a_1=0,\ b_s=n+1,\
a_i<a_{i+1},\ b_i<b_{i+1},\
b_i-a_i\ge k,\
0\le b_i-a_{i+1}<k-1
\]
(here $0$ and $n+1$ are formal endpoints).
Impose the interval rank conditions
\[
\begin{aligned}
r_V\bigl([a_i,b_i]\cap[1,n]\bigr)&\le k-1
&&\text{for every block},\\
r_V([a_{i+1},b_i])&\le b_i-a_{i+1}
&&\text{for every overlap}.
\end{aligned}
\]
Every consecutive $k$-column interval either lies in a block or contains an overlap.
It is shown in \cite{HS04} that these sequences give all the irreducible components of 
the consecutive-minor locus, without repetitions. After restricting to rank-$k$ matrices, 
each locus is an interval positroid variety: its generic ranks on all intervals specify 
its interval rank matrix. The one-block sequence $\Gamma=([0,n+1])$ requires 
$r_V([1,n])\le k-1$, so it has no points in $\mathrm{Gr}_{k,n}$. Every other sequence 
gives one distinct component of ${\mathsf Y}_{k,n}$.

Finally, \cite{HS04} proves that the number of such $\Gamma$ is equal to $F^{(k)}_{n-k+1}$
(using the recursive description).
\end{proof}

\begin{example}
Let $(k,n)=(3,5)$. The special fiber
${\mathsf Y}_{3,5}\subset\Gr_{3,5}$ is defined by
$P_{123}=P_{234}=P_{345}=0$.
The compositions of $n-k+1=3$ with parts at
most $3$, other than $(1,1,1)$, are
$(3)$, $(1,2)$, $(2,1)$.
The three corresponding interval positroid
varieties have rank sets
\[
\begin{aligned}
\mathcal R_1=\{[1,5],[3,3]\},\quad \mathcal R_2=\{[1,3],[3,4]\},\quad 
\mathcal R_3=\{[2,3],[3,5]\}.
\end{aligned}
\]

Writing a point of $\Gr_{3,5}$ as a matrix
with columns $v_1,\ldots,v_5$, these
components can equivalently be described by
\[
\begin{aligned}
X_1
&=\{V:v_3=0\},\\
X_2
&=\{V:
r_{1,3}(V)\le 2,\ 
r_{3,4}(V)\le 1,\
r_{1,4}(V)\le 2\},\\
X_3
&=\{V:
r_{2,3}(V)\le 1,\
r_{3,5}(V)\le 2,\
r_{2,5}(V)\le 2\},
\end{aligned}
\]
where
$r_{i,j}(V)=\operatorname{rank}(v_i,\ldots,v_j)$.
Each component has dimension $3$, and
${\mathsf Y}_{3,5}=X_1\cup X_2\cup X_3$.
Thus the number of irreducible components is
$F_3^{(3)}-1=4-1=3$.
\end{example}

\begin{cor}
All the maximal elements of the poset of positroids not containing $C_a$, $a\in [n-k+1]$ are 
interval positroids of codimension $n-k+1$.  
\end{cor}

\end{document}